\documentclass[12pt,reqno]{amsart}

\usepackage{amsmath,amsthm,amsfonts,amscd,flafter,epsf}

\newcommand{\stdspace}{\hskip 0.35em plus 0.15em \ignorespaces}
\let\qua\stdspace

\newtheorem{theorem}{Theorem}[section]

\newtheorem{conj}[theorem]{Conjecture}

\newtheorem{remark}[theorem]{Remark}

\newcommand\Rk{\mathrm{rk}}

\newcommand\os{Ozsv\'ath-Szab\'o }

\newcommand{\C}{\mathbb{C}}

\newcommand{\Z}{\mathbb{Z}}

\newcommand{\cm}{\cdot}

\newcommand{\ModSWfour}{\mathcal{M}}
\newcommand{\ModFlow}{\ModSWfour}

\newcommand{\SpinC}{{\mathrm{Spin}}^c}

\newcommand\Hom{\mathrm{Hom}}

\newcommand\abuts\Rightarrow
\newcommand\Sym{\mathrm{Sym}}

\newcommand\ModSphere{\ModFlow\left({\mathbb S}\longrightarrow
\Sym^{g-1}(\Sigma_{1})\times \Sym^2(\Sigma_{2})\right)}
\newcommand\ModSpheres\ModSphere

\newcommand\HFp{\HFb}

\newcommand\HFa{\widehat{HF}}
\newcommand\HFb{HF^+}

\newcommand\UnparModSp{\widehat \ModSp}
\newcommand\UnparModFlow\UnparModSp
\newcommand\Mod\ModSp

\newcommand\PD{\mathrm{PD}}

\newcommand{\spinct}{\mathfrak t}

\newcommand\ModMaps{\mathcal M}
\newcommand\ModSp\ModMaps

\newcommand\del{\partial}

\newcommand\Dual{\mathcal D}
\newcommand\Duality\Dual

\newcommand\Char{\mathrm{Char}}
\newcommand\Combp{{\mathbb H}^+}
\newcommand\DCombp{{\mathbb K}^+}

\newcommand\NatTransp{T^+}
\newcommand\InjMod[1]{{\mathcal T}^+_{#1}}

\begin{document}

\title{On Plumbed $L$-spaces}
\author{Raif Rustamov}
\address{The Program in Applied and Computational Mathematics, Princeton University\\New Jersey 08540, USA}
\email{rustamov@princeton.edu}
\begin{abstract}
In the quest for $L$-spaces we consider links of isolated complete
intersection surface singularities. We show that if such a manifold
is an $L$-space, then it is a link of a rational singularity. We
also prove that if it is not an $L$-space then it admits a
symplectic filling with $b_2^+>0$. Based on these results we pin
down all integral homology sphere $L$-spaces in this realm.
\end{abstract}
\maketitle

\section{Introduction}
To understand the geometric content of the \os Floer homology
better, it is important to investigate the class of three-manifolds
with Floer homology of a particularly simple form:
\emph{$L$-spaces}. A rational homology sphere $Y$ is called an
$L$-space if $\Rk (\HFa(Y)) = |H_1(Y)|$. Every Lens space is an
$L$-space - this is where the name comes from. However, Lens spaces
are not the only representatives of this class, rather there is a
very big supply of $L$-spaces. Indeed, for a knot $K$ in $S^3$ such
that $S^3_p(K)$ is an $L$-space with $p>0$, $S^3_{p+1}(K)$ is also
an $L$-space. For example, letting $K$ to be the pretzel knot
$P(-2,3,7)$ one sees that $S^3_{18}(K)$ is an $L$-space. It follows
that for any integer $n \geq 18$, $S^3_n(K)$ is an $L$-space. It is
worth mentioning that the knot $P(-2,3,7)$ is hyperbolic.

In the quest for $L$-spaces let us consider links of isolated
surface singularities. By the work of N\'emethi \cite{Nemethi} it is
known that links of rational surface singularities are $L$-spaces.
We prove the converse of this statement for the singularities that
are complete intersections.

\begin{theorem}
\label{ll} Let $Y$ be a link of isolated complete intersection
surface singularity. If $Y$ is an $L$-space, then it is the link of
a rational surface singularity.
\end{theorem}
Note that none of the $L$-spaces constructed as described in the
beginning is an integral homology sphere. In fact, we know only a
limited number of such $L$-spaces, and one is lead to the following
conjecture.

\begin{conj} (Ozsv\'ath and Szab\'o) If an irreducible integral homology
sphere $Y$ is an $L$-space, then
 $Y=S^3$ or
$Y=\pm \Sigma(2,3,5)$.
\end{conj}

An irreducible manifold is one that cannot be written as a connected
sum of two manifolds different from $S^3$. This is why the integral
homology sphere $L$-spaces $$n(\Sigma(2,3,5))\#m(-\Sigma(2,3,5)),$$
where $n, m \geq 0$, are not included in the list. We provide more
evidence for this conjecture by proving the following theorem.

\begin{theorem} \label{aa}
Consider an irreducible manifold $Y \neq S^3$ which is a link of an
isolated complete intersection surface singularity. If $Y$ is an
$L$-space, then $\pm Y=\Sigma(2,3,5)$.
\end{theorem}

The following theorem from \cite{Genus} is a good indication of how
much Heegaard Floer homology, namely being an $L$-space, affects the
geometry of a manifold in general.

\begin{theorem}
\label{foil} (Ozsv\'ath and Szab\'o) All symplectic fillings of an
$L$-space have $b^+_2=0$.
\end{theorem}

Restricting to complete intersection surface singularities we can
prove the converse statement.

\begin{theorem}\label{sympfill}
If a link of isolated complete intersection surface singularity is
not an $L$-space, then it bounds a Stein four-manifold with
$b_2^+>0$. In particular, it has a symplectic filling with
$b_2^+>0$.
\end{theorem}

The paper is organized as follows: We start with the review of the
plumbing construction and its relation to surface singularities.
Section 4 contains the proofs of theorems 1.1 and 1.5. The last
section starts with the review of the combinatorial method for
computing the \os Floer homology in the case of almost rational
graphs and concludes with the proof of 1.3.

\medskip
\noindent {\bf Acknowledgments} \qua I am pleased to thank my
advisor Zolt\'an Szab\'o whose kind attitude and outstanding
teaching made this work very enjoyable. I would like to express my
gratitude to Andr\'as N\'emethi who spotted an error in the first
version of this paper and gave many suggestions to improve the
paper.

\section{The Plumbing Construction}
 Let $G$ be a weighed forest, i.e. not necessarily
connected graph without any cycles, to each vertex of which an
integer is assigned. In what follows, the term "graph" will be
used for this specific class of objects. Let $m(v)$ and $d(v)$ be
respectively the weight and degree of the vertex $v$. To each
vertex $v$ associate the disk bundle on $S^2$ with Euler number
equal to $m(v)$. For any edge $vw$ of the graph, choose a small
disk on the spheres corresponding to $v$ and $w$. When the bundles
are restricted to these disks, (topologically) we have $D^2 \times
D^2$. Glue the two bundles along this two $D^2 \times D^2$'s, by
exchanging the base and fiber coordinates -- i.e. \emph{plumb}
these two bundles. We denote the resulting four-manifold with
boundary by $X(G)$. It is said to be obtained by plumbing on $G$.
We let $Y(G)$ be the oriented boundary of $X(G)$. Three-manifolds
obtained so are called a \emph{plumbed} manifolds.

Note that $X=X(G)$ is simply-connected, and has a very natural
representation for its second homology group furnished by the zero
sections of the bundles. In fact, the group $H_2(X;\Z)$ is the
lattice freely spanned by vertices of $G$. Denoting by $[v]$ the
homology class in $H_2(X;\Z)$ corresponding to a vertex $v$ of
$G$, the values of the intersection form of $X$ on the basis are
given by $[v]\cm[v] = m(v)$; $[v]\cm[w] = 1$ if $vw$ is an edge of
$G$ and $[v]\cm[w]=0$ otherwise. $G$ is called
\emph{negative-definite} if the form is negative-definite.

The plumbed manifold $Y(G)$ has an alternative description in
terms of surgery on links. To each vertex $v$  of $G$ one
associates an unknot in $S^3$, framed by $m(v)$. Two unknots are
linked geometrically once if there is an edge between the
corresponding vertices. Since our graph does not have any cycles,
this can be done unambiguously. Making surgery on this link gives
us $Y(G)$. The four-manifold $X(G)$ can be obtained by attaching
2-handles to the 4-ball according to this link. It is clear that
we can describe the same manifolds by different graphs, and Kirby
moves would  allow to move us between any two plumbing
representations. These Kirby moves can be conveniently described
as some combinatorial changes of the plumbing graph, see for
example \cite{moves}. A graph $G$ will be called minimal if it
cannot be made smaller (less vertices and/or edges) using just one
Kirby move.
\section{Plumbed manifolds and links of singularities} Given $m$
analytic functions $f_1, ..., f_m:(\C^n,0) \rightarrow (\C,0)$ one
can consider their common zero locus $Z=\{f_1=0,...,f_m=0\}$. For
each point $p \in Z$ consider the matrix $$J(p) = \left(\frac{\del
f_i}{\del z_j} \right) (p).$$ We will be interested in the case
when the rank of $J(p)$ is equal to $n-2$ for every $p \in
Z-\{0\}$. If the rank of $J(0)$ is also equal to $n-2$, $Z$ is
said to be smooth and is analytically isomorphic to $(\C^2,0)$.
Otherwise, $(Z,0)$ is called an \emph{isolated surface
singularity}.

Assuming that $(Z,0)$ is irreducible and normal, one can associate
to it a three-manifold as follows. Given an embedding of $(Z,0)$
into $(\C^n,0)$, one can find $\epsilon_0$ such that for any $0 <
\epsilon \leq \epsilon_0$, the $2n-1$ dimensional sphere of radius
$\epsilon$ centered at the origin intersects $(Z,0)$
transversally. This intersection is an oriented three-manifold
independent of the arbitrary choices. It is called the \emph{the
link of singularity}  of $(Z,0)$. If $(Z,0)$ is smooth the link is
$S^3$.

Fix a \emph{good resolution} $\pi : (\tilde{Z},E) \rightarrow
(Z,0)$. Then the exceptional set $E=\pi^{-1}(0)$ is a union of
smooth curves intersecting transversally. One can consider the
\emph{resolution dual graph} $\Gamma$ of $E$, and it turns out
that the link of singularity of $(Z,0)$ can be identified with the
 with the plumbed three-manifold $Y(\Gamma)$. The
link is a rational homology sphere if and only if the components
of $E$ are rational, and the graph $\Gamma$ does not have any
loops. As a result, all links of singularities which are rational
homology spheres are plumbed manifolds. Conversely, a plumbed
manifold can be realized as a link of some normal surface
singularity if and only if it can be obtained as a plumbing on a
negative definite tree.

The \emph{geometric genus} $p_g(Z,0)$ of a singularity is defined
as the dimension of $H^1(\tilde{Z}, \mathcal{O}_{\tilde{Z}})$. If
$p_g=0$, the singularity is called \emph{rational}. If a negative
definite plumbed manifold is the link of a rational singularity,
then the plumbed manifold and the plumbing graph are said to be
rational. A graph $G$ is called \emph{almost-rational} if  by
decreasing just one of the weighs, the graph can be made rational,
see \cite{Nemethi}.

\section{Proofs of theorems \ref{ll} and \ref{sympfill}}
Recall that a smooth compact oriented $2n$-dimensional manifold with
possibly non-empty boundary is said to be \emph{symplectic} if it
admits a closed 2-form $\omega$ such that $\omega^n$ does not
vanish. A complex manifold $M$ (possibly with boundary) is a
\emph{Stein} manifold if it admits a positive, proper, and strictly
plurissubharmonic (Morse) function for which $M$ is a regular level
set. By a theorem of Grauert, a complex manifold is Stein if and
only if it embeds analytically and properly into $\C^N$ for some
$N$. Note that Stein manifolds are symplectic.

A 1-form $\alpha$ on an oriented closed three-manifold $Y$ is called
a \emph{contact form} if $\alpha \wedge d\alpha>0$. Given a contact
manifold , i.e. a pair $(Y,\alpha)$, it said to admit a symplectic
filling if there is a symplectic four-manifold $(X, \omega)$ such
that $\del X=Y$  and $\omega|_{\mathrm{ker} \alpha} \neq 0$.

Let $Y$ be  the link of an isolated complete intersection surface
singularity, say $(Z,0)$, given by $f=(f_1,...,f_m)$. To complete
intersection surface singularities one can associate a four-manifold
$F$ called \emph{Milnor fiber}. It is defined as the manifold given
by intersection $f^{-1}(\delta) \cap B(\epsilon)$, where
$B(\epsilon)$ is a sufficiently small ball about the origin,
$\delta$ is a general point of $C^m$ very close to $0$. $F$ is a
smooth simply-connected four-manifold. It is oriented by the complex
structure on its interior. Furthermore, its boundary is
diffeomorphic to the link of singularity of $(Z,0)$, i.e. $Y$. Note
that $F$ is Stein, because it is a piece of complex surface in a
ball.

There is a connection between the geometric genus of singularity
$(Z,0)$ and the Betti number of its Milnor fiber $F$. In fact, we
have the following formula due to Wahl, Durfee and Steenbrink:
$$p_g(Z,0)=b_2^+(F)/2.$$
The proof of this formula in its whole generality uses mixed Hodge
theory.

We will use the following theorem.
\begin{theorem}{(\cite{Genus})}
\label{stein} Let $Y$ be an $L$-space. If $X$ is any Stein
four-manifold with boundary $\pm Y$, then $b^+_2(X)=0$.
\end{theorem}

To prove the theorem \ref{ll}, suppose that $Y$ is an $L$-space.
Consider the Milnor fiber $F$ of the corresponding singularity
$(Z,0)$. Since $F$  is a Stein manifold and $\del F=Y$, we have that
$b_2^+(F)=0$. By the formula of Wahl, Durfee and Steenbrink we get
$p_g(Z,0)=b_2^+(F)/2=0$, thus $(Z,0)$ is a rational singularity.
This concludes the proof of the statement that if a link of a
complete intersection surface singularity is an $L$-space, then is
is a link of rational singularity.

Now we turn to the theorem \ref{sympfill}. If $Y$ is not an
$L$-space then by the previous argument $(Z,0)$ is not rational. As
a result $b^+_2(F)=2p_g(Z,0)>0$. By a result of Varchenko
\cite{Varchenko}, $Y(G)$ admits a special contact structure $\alpha$
which can be shown to be the contact boundary of the Milnor fiber.
Thus $F$ provides us with a Stein filling of $(Y, \alpha)$ with
$b^+_2>0$. This finishes the proof.

\section{Integral homology spheres among links of complete
intersection surface singularities}

\subsection{Heegaard Floer Homology of Almost Rational Plumbings}
\label{alg} Here we review the results of \cite{Plumbing} and
\cite{Nemethi}. Let $\InjMod{0}$ denote the $\Z[U]$-module which is
the quotient of $\Z[U,U^{-1}]$ by the submodule $U\cm \Z[U]$, graded
so that the element $U^{-d}$ (for $d\geq 0$) has degree $2d$. Let
$\InjMod{d}=\InjMod{0}[d]$, i.e. the grading is shifted by $d$.

Denoting by $\Char(G)$  the set of characteristic vectors for the
intersection form of $X(G)$, define
$$\Combp(G)\subset \Hom(\Char(G),\InjMod{0})$$
to be the set of finitely supported functions satisfying the
following relations for all characteristic vectors $K$ and vertices
$v$:

\begin{equation}
\label{eq:AdjRel} U^n\cm \phi(K+2\PD[v]) = \phi(K),
\end{equation}
if $2n=\langle K,v \rangle + v\cm v \geq 0$; and
\begin{equation}
\label{eq:AdjRel2} \phi(K+2\PD[v]) = U^{-n}\cm \phi(K)
\end{equation}
for $n<0$.

We can decompose $\Combp(G)$ according to $\SpinC$ structures over
$Y$. Note that the first Chern class gives an identification of the
set of $\SpinC$ structures over $X=X(G)$ with the set of
characteristic vectors $\Char(G)$. Observe that the image of
$H^2(X,\partial X;\Z)$ in $H^2(X;\Z)$ is spanned by the Poincar\'e
duals of the spheres corresponding to the vertices. Using the
restriction to boundary, it is easy to see that the set of $\SpinC$
structures over $Y$ is identified with the set of $2 H^2(X,\partial
X;\Z)$-orbits in $\Char(G)$.

Fix a $\SpinC$ structure $\spinct$ over $Y$. Let $\Char_\spinct(G)$
denote the set of characteristic vectors for $X$ that can be
realized as the first Chern classes of $\SpinC$ structures whose
restriction to the boundary is $\spinct$. Similarly, we let
$$\Combp(G,\spinct)\subset \Combp(G)$$ be the subset of maps
supported on $\Char_\spinct(G)\subset \Char(G)$. We have a direct
sum splitting: $$\Combp(G)\cong \bigoplus_{\spinct\in\SpinC(Y)}
\Combp(G,\spinct).$$

The grading on $\Combp(G)$ is introduced as follows: we say that an
element $\phi\in \Combp(G)$ is homogeneous of degree $d$, if for
each characteristic vector $K$ with $\phi(K)\neq 0$, $\phi(K)\in
\InjMod{0}$ is a homogeneous element with:
\begin{equation}
\label{eq:DefOfDegree} \deg(\phi(K))-\frac{K^2+|G|}{4}=d.
\end{equation}

The following theorem is proved in \cite{Nemethi}. It is a
generalization of a result from \cite{Plumbing}.
\begin{theorem}
\label{cucuk} Let $G$ be an almost rational negative-definite
weighted forest. Then, for each $\SpinC$ structure $\spinct$ over
$-Y(G)$, there is an isomorphism of graded $\Z[U]$-modules,
$$\HFp(-Y(G),\spinct)\cong \Combp(G,\spinct).$$
\end{theorem}

For calculational purposes it is helpful to adopt the dual point of
view. Let $\DCombp(G)$ be the set of equivalence classes of elements
of $\Z^{\geq 0} \times \Char(G)$ (and we write $U^m\otimes K$ for
the pair $(m, K)$) under the following equivalence relation. For any
vertex $v$, let
$$2n=\langle K,v \rangle + v\cm v.$$
If $n\geq 0$, then
\begin{equation}
\label{equation:rel1} U^{n+m}\otimes (K+2\PD[v]) \sim U^m\otimes K,
\end{equation}
while if $n\leq 0$, then
\begin{equation}
\label{equation:rel2} U^m\otimes (K+2\PD[v]) \sim U^{m-n}\otimes K.
\end{equation}

Starting with a map $$\phi\colon \Char(G)\longrightarrow
\InjMod{0},$$ consider an induced map $${\widetilde \phi}\colon
\Z^{\geq 0}\times \Char(G)\longrightarrow \InjMod{0}$$ defined by
$${\widetilde \phi}(U^n\otimes K)=U^n\cm \phi(K).$$ Clearly, the
set of finitely-supported functions $\phi\colon
\Char(G)\longrightarrow \InjMod{0}$ whose induced map ${\widetilde
\phi}$ descends to $\DCombp(G)$ is precisely $\Combp(G)$.

A \emph{basic element} of $\DCombp(G)$ is one whose equivalence
class does not contain  any element of the form $U^m\otimes K$ with
$m>0$. Given two non-equivalent basic elements $K_1 = U^0 \otimes
K_1$ and $K_2 = U^0 \otimes K_2$ in the same $\SpinC$ structure, one
can find positive
 integers $n$ and $m$ such that
$$ U^{n} \otimes K_1 \sim U^{m} \otimes K_2.$$
If, moreover, the numbers $n$ and $m$ are as small as possible, then
this relation will be called the \emph{minimal relationship} between
$K_1$ and $K_2$. One can see that $\DCombp(G)$ is specified as soon
as one finds its basic elements and the minimal relationships
between each pair of them.

\begin{remark}The manifold under consideration
is an $L$-space if and only if the number of basic vectors is equal
to the number of the $\SpinC$ structures. Of course, in this case
there will be one basic vectors in each $\SpinC$ structure, and the
basic vectors will be unrelated.
\end{remark}

Now we review the algorithm given in \cite{Plumbing} for calculating
the basic elements.

Let $K$ satisfy
\begin{equation}
\label{eq:PartBox} m(v)+2 \leq \langle K, v\rangle \leq -m(v).
\end{equation}

Construct a sequence of vectors $K=K_0,K_1,\ldots,K_n$, where
$K_{i+1}$ is obtained from $K_i$ by choosing any vertex $v_{i+1}$
with
$$\langle K_i,v_{i+1}\rangle = -m(v_{i+1}),$$
and then letting
$$K_{i+1}=K_i+2\PD[v_{i+1}].$$
Note that any two vectors in this sequence are equivalent.

This sequence can terminate in one of two ways: either
\begin{itemize}
\item the final vector $L=K_n$ satisfies the inequality,
\begin{equation}
\label{eq:OtherPartBox} m(v) \leq \langle L, v\rangle \leq -m(v)-2
\end{equation}
at each vertex $v$ or \item there is some vertex $v$ for which
\begin{equation}
\label{eq:OtherTermination} \langle K_{n},v \rangle > -m(v).
\end{equation}
\end{itemize}

It turns out that the basic elements of $\DCombp(G)$ are in
one-to-one correspondence with initial vectors $K$ satisfying
inequality~\eqref{eq:PartBox} for which the algorithm above
terminates in a characteristic vector $L$ satisfying
inequality~\eqref{eq:OtherPartBox}.

\subsection{$L$-spaces and rational singularities} Nemethi proves in
\cite{Nemethi} the following theorem.

\begin{theorem}{(N\'emethi)}
For a negative definite plumbing graph $G$ the following are
equivalent:

\noindent (i) $G$ is rational;

\noindent(ii) $\Combp(G,\spinct_{can})=\InjMod{d}$ for some $d$;

\noindent(iii) For every $\spinct$, $\Combp(G,\spinct)=\InjMod{d}$
for some $d=d(\spinct)$.
\end{theorem}

Here $\spinct_{can}$ is the \emph{canonical $\SpinC$ structure} on
$Y(G)$. In our description it can be specified as follows. Let
$K_{can}$ be the \emph{canonical characteristic vector}, i.e. the
vector in $\Char(G)$ that satisfies $$\langle K, v\rangle =
m(v)+2,$$ for every vertex $v$ of the plumbing graph $G$. There is a
unique $\SpinC$ structure, which we will denote by $\spinct_{can}$,
such that the inclusion $K_{can} \in \Char(G,\spinct_{can})$ is
true.

Recall that for each $\SpinC$ structure $\spinct$, we have the map
of $\Z[U]$-modules $$\NatTransp\colon \HFp(-Y(G),\spinct)
\longrightarrow \Combp(G,\spinct).$$ N\'emethi proves that if $G$ is
rational then $\NatTransp$ is an isomorphism. This proves the
following theorem.

\begin{theorem}{(N\'emethi)} For any  negative definite
\emph{rational} plumbing graph $G$, the three-manifold $Y(G)$ is an
$L$-space, i.e. links of rational surface singularities are
$L$-spaces.
\end{theorem}

\subsection{Proof of Theorem \ref{aa}}
By the theorem \ref{ll} we have to find all rational graphs for
which the plumbed manifold is an integral homology sphere. Note that
by the condition of irreducibility we can concentrate on connected
graphs, because for disconnected graphs the resulting manifold is
the connected sum of plumbed manifolds corresponding to the
components of the graph. Without loss of generality we assume that
the graph to be minimal. We claim that the only connected rational
negative definite plumbing graph which gives rise to an integral
homology sphere is (negative) $E_8$. Obviously, this claim implies
the theorem.

We will use the characterization of rational graphs by Nemethi, and
the fact that any subgraph of a rational graph is also rational. The
same is true about being negative definite. Letting $G$ be our
weighed plumbing graph with weigh at the vertex $v$ denoted by
$m(v)$, we have several possibilities.

\medskip \noindent {\bf{Case 1}} \qua The graph has at least one
vertex, say $w$, with weigh $-1$. The condition of minimality and
the fact that we are not interested in $S^3$ means that there are at
least three vertices adjacent to $w$. Consider the graph $G'$ made
of the vertex $w$ and all vertices adjacent to it. We will arrive at
a contradiction by assuming that $G$ is rational. Note that by this
assumption, $G'$ is also rational. Thus $\pm Y(G')$ is an $L$-space.
Consider the graph $-G'$, i.e. the same graph with weighs multiplied
by $(-1)$. Using usual transformations that do not change the
plumbed manifold, we can modify this graph to get a new graph $H$
with all of its weighs $\leq -2$. $H$ will have a strand of $-2$'s
of length $m(v)-1$ for each vertex $v\neq w$. These strands will
join together at a central vertex with weigh
$m=1-\mathrm{\#}(\mathrm{vertices \stdspace adjacent \stdspace to
\stdspace} w \mathrm{\stdspace in \stdspace} G')$. Thus we have
$Y(G')=-Y(H)$. Consider the four-manifold $X(H)$. It has $b^+_2=1$
and it is Stein, because it can be obtained by Legandrian surgery on
a link in $S^3$. The last is true because all the weighs of $H$ are
less than $-1$. Remembering that $\del X(H) = Y(H)$, the theorem
\ref{stein} implies that $\pm Y(H)=\mp Y(G')$  is not an $L$-space,
thus yielding a contradiction.

\medskip \noindent {\bf{Case 2}} \qua The graph has no vertex with weigh
$-1$, and there is at least one vertex $w$ with $m(w)\leq -3$.
Consider the canonical vector $K_{can}$. We can modify it by adding
$2$ to its entry corresponding to the vertex $w$. Let $K'$ denote
this new vector. Obviously $K'$ satisifies the inequality
\ref{eq:PartBox}. Now if we run the algorithm described in the
section \ref{alg} on this vector, the algorithm will stop at the
very first step yielding $L=K'$. Clearly $L$ satisfies the
inequality \ref{eq:OtherPartBox}, thus $K'$ is a basic vector. Since
we want $G$ to be rational, it follows that $K'$ cannot be in the
same $\SpinC$ structure as $K_{can}$, because otherwise for no $d$,
the identification $\Combp(G,\spinct_{can}) = \InjMod{d}$ would be
possible. As a result, we have at least two different $\SpinC$
structures on $\pm Y(G)$, which means that $\pm Y(G)$ is not an
integral homology sphere.

\medskip \noindent{\bf{Case 3}} \qua The only remaining case to consider is
when all of the weighs on the graph are $-2$. It is a simple
exercise in combinatorics to see that the condition of negative
definiteness leaves us only one choice -- the (negative) $E_8$
plumbing. Indeed, no such graph can have a vertex with degree
greater than 3. Neither it can contain two vertices with degrees
both greater than 2. The remaining case of the star graph with three
emanating strings is easily investigated.

\end{document}